\documentclass[a4article,12pt]{article}

\usepackage[T2A]{fontenc}
\usepackage[russian,english]{babel}
\usepackage[cp1251]{inputenc}
\usepackage{amsthm,amsmath}

\usepackage{amssymb}
\usepackage{graphicx}
\usepackage{color}

\makeatletter
\renewcommand{\@biblabel}[1]{#1.}
\makeatother

\begin {document}
%-------------------------------------------------------------------------

%
%

\begin{center}
\textbf{Spectral analysis of one-term symmetric differential operators of even order with interior singularity}
\end{center}

\centerline{Irina Braeutigam \footnote{The author acknowledges the support of the Ministry of Education and Science
of Russian Federation and the German Academic Exchange Service (DAAD) under the program
“Mikhail Lomonosov” (grant no. 325-A/13/74976) and the Russian Science Support Foundation.}}

\bigskip

\textbf{Abstract.}
In this paper we discuss the spectral properties of one-term symmetric differential operators of even order with interior singularity, namely, we determine the deficiency numbers, describe its self-adjoint extensions and their spectrum. It is assumed that the operators are generated by the differential expression

$l_{2m}[y](x)=(-1)^m(c(x)y^{(m)})^{(m)}(x)$, where $x \in I:=[-1,1]$,

\noindent the coefficient $c(x)$ has one zero on the set $I$ and the orders of this zero on the right side and the left side are not necessarily equal.



\numberwithin{equation}{section}
\newtheorem{theorem}{Theorem}[section]
\newtheorem{lemma}[theorem]{Lemma}

\allowdisplaybreaks

\section{Introduction}
Differential operators with interior singularity occur in physical applicationы and many areas of contemporary analysis. Although they were mentioned already in the book \cite{DanSchw}, spectral analysis of such operators is surprisingly rarely examined. One of the first works in this direction was the paper of J.P.Boyd \cite{Boyd}, where he considered Sturm-Liouville operators with an interior pole. Later W.N. Everitt and A.Zettl \cite{EvZet1} developed a theory of self-adjoint realization of Sturm-Liouville problems on two intervals in direct sum of Hilbert spaces associated with these intervals. In \cite{EvZet2} they extended this theory to higher order differential operators and any number of finite or infinite intervals. In the last few years some articles appeared  where  self-adjoint domains of ordinary differential operators in direct sum of Hilbert spaces were described in terms of real-parameter solutions of the corresponding differential equations, see, for example, \cite{Sun},\cite{Wang}. In the papers mentioned above, minimal and maximal differential operators, generated by differential expressions on more than one interval are described only in the sense of the direct sum of operators given on each (sub)interval. However, if we consider adjacent intervals then there is also another approach to describe the minimal differential operators which began to develop Y.B. Orochko. Namely, this approach is based on the next heuristic interpretation. The interior singularity could be considered as an interior barrier for two different evolution processes. If these processes can not propagate to the adjacent interval over this point than the minimal differential operator, generated on the whole interval has decomposition into an orthogonal direct sum of the minimal operators, generated on each subinterval. In other cases such decomposition is not possible and the minimal operator on the whole interval is a symmetric extension of the orthogonal direct sum of the minimal operators, generated on each subinterval. More details could be found in \cite{Oro3}. In the present work we will use this approach to describe one-term minimal differential operators with interior singularity.

Let us consider an ordinary differential expression of the arbitrary order $2m$ ($m=1,2, \ldots$)
\begin {equation}
\label {osn}
l_{2m}[y](x)=(-1)^m(c(x)y^{(m)})^{(m)}(x), \quad x \in I:=[-1,1].
\end{equation}

 We introduce now the concept of an irregular differential expression and an extended concept of the order of the zero of the coefficient $c(x)$.

Analogous to the book~\cite[Ch. XIII]{DanSchw},
we call the differential expression  $l_{2m}[y]$  irregular
differential expression when the  coefficient $c(x)$ of this expression vanishes at some points of set I.

Point $x_0 \in I$ is called  right zero of the coefficient  $c(x)$ of the order $p>0$
(and respectively left zero of the order $q>0$), if $c(x)=(x-x_0)^pa(x)$ for $x \in (x_0;x_0+h]$, $h>0$, where
$a(x)$ is a positive or negative function on the segment $[x_0;x_0+h]$
(respectively $c(x)=|x-x_0|^qb(x)$ for $x \in [x_0-h;x_0)$, where $b(x)$ is a positive or negative  function on $[x_0-h;x_0]$).

  Suppose that the coefficient $c(x)$ of the expression $l_{2m}$ is determined on $I$ and has on this set a single zero $x_0=0$ of the right order $p$ and the left order
$q$, where $p,q \in \{1,2, \ldots, 2m-1\}$, which means it can be represented as
$$
c(x)= \begin{cases}
x^p a(x), &\text { if $x \in [0,1]$,}\\
(-x)^q b(x) , &\text { if $x \in [-1,0]$,}
\end{cases}
$$
  while the functions $a(x), b(x)$ are real-valued functions on $I$ and can be represented as power series which are convergent whenever  $|z|< 1$
$$
 a(z):=a_0 + \sum\limits_{j=1}^{+\infty} a_jz^j , \;\; a_0 \ne 0,
  $$
$$
 b(z):=b_0 + \sum\limits_{j=1}^{+\infty} b_jz^j , \;\; b_0 \ne 0.
  $$
Here we note that the further results given in this paper also remain true when the functions  $a(z), b(z)$ are analytical with $|z| \le x_0 <1$ for an  $x_0 \in (0,1]$, Lebesgue integrable outside $[-x_0;x_0]$ and do not vanish there. The case $x_0=1$ is chosen for convenience.

The differential expression (\ref{osn}) has one interior singularity $x=0$ due to the fact that $\frac  1 {c(x)}$ is not Lebesgue integrable in left and right neighborhoods of this point. Therefore we also can call (\ref{osn}) the differential expression with interior singularity.

We define the quasi-derivatives
of a given function $f(x)$ which correspond to the expression $l_{2m}$ as follows:
$$
f^{[l]}(x)= \begin{cases}
f^{(l)}(x), &\text { if $l=0,1, \ldots m-1$,}\\
(-1)^m(c(x)f^{(m)})^{(l-m)}(x) , &\text { if $l=m, m+1 \ldots 2m-1$.}
\end {cases}
$$
Let  $f$ and $g$ be two functions for which the expression $l_{2m}$ is defined. As it is well known, Green's formula applies  (~\cite[ch.V, \S 15]{Naim})
$$
\int\limits_\alpha^\beta (l_{2m}[f] \overline g - f l_{2m} [\overline {g}])dx = [f,g]|^\beta_\alpha, \qquad \alpha,\beta \in I,
$$
where
\begin {equation}
\label {biform}
[f,g](x)= (-1)^m\sum\limits_{j=1}^{2m}(-1)^{j-1} \{ f^{[2m-j]}(x)\overline
g^{[j-1]}(x)\}=(-1)^mG^*EF
\end {equation}
and $$
G=\begin{pmatrix}
g \\
g^{[1]} \\
\vdots \\
g^{[2m-1]}
\end{pmatrix}, \qquad
F=\begin{pmatrix}
f \\
f^{[1]} \\
\vdots \\
f^{[2m-1]}
\end{pmatrix}, \qquad E=((-1)^{r-1} \delta_{r,2m+1-s})_{r,s=1}^{2m}.
$$
Let us note here that the sesquilinear form
$[f,g](x)$ is defined by formula~(\ref {biform}) for $x \in [-1,0) \cup (0,1]$ and  for any two functions $f$ and $g$ which have absolutely continuous quasi-derivatives of all orders up to and including the order of $2m-1$ in the neighborhood of point $x$, which is a part of the set $[-1,0) \cup (0,1].$ We denote
 $$
[f,g](-0)=\lim_{x \to -0} [f,g](x), \;\; [f,g](+0)=\lim_{x \to
+0}[f,g](x),
 $$
provided that the indicated limits exist.

Next we define a minimal closed symmetric differential operator in $L_2(I)$.

We denote by ${D}^{\prime}_0$ the set of all  infinitely differentiable functions  $y$ on $I$ which vanish identically outside a finite set $[\alpha, \beta]  \subset [-1,1]$; this set may be different for different functions.  By $L^{\prime}_0$ we denote the operator with the domain ${D}^{\prime}_0$ . So, for $y \in {D}^{\prime}_0$ we have $L^{\prime}_0y=l_{2m}[y]$.  According to Green's formula, the operator  $L^{\prime}_0$ is a symmetric operator with an  everywhere dense domain
${D}^{\prime}_0$
and consequently permits closure in the space $L_2(I)$.

We denote the closure of $L^{\prime}_0$ as $L_0^{pq}$. Thus, $L_0^{pq}$ is a minimal closed symmetric operator generated by the irregular differential
expression $l_{2m}$ in $L_2(I)$.
Let us  denote its domain with the symbol ${D}_0$.

There is another equivalent characterization of $D_0$ in terms of the sesquilinear forms (\ref{biform}), namely,
$$D_0=\{ f\in D | f^{[k]}(1)=f^{[k]}(-1)=0, k=0, \ldots, 2m-1,$$
$$[f,g](0-)-[f,g](0+)=0, \forall g \in D \}.$$

Since $L_0^{pq}$ is a real operator, its deficiency numbers in the upper and lower open complex semiplanes are equal. We denote their common value by $n_{pq}$.

It is well known that if $c(x)>0$ or $c(x)<0$ at a certain interval, the deficiency numbers of the minimal closed symmetric operator do not exceed the order of the corresponding differential expression. In~\cite{Oro4} Y.B. Orochko  gives examples showing that if an interval has a finite or countable set of zeros of the  coefficient $c(x)$, then the opposite inequality can be true.

Furthermore, in ~\cite{Oro5} he considered the minimal symmetric differential operator $L_0$ generated by the irregular differential
expression~(\ref {osn}) on $I$,
 $c(x):=x^pa(x)$, $p \in \{1,2, \ldots, 2m-1\}$, $a(x)$ is an infinitely differentiable real-valued function and $a(x) \ne 0$ for any
$x \in I$. It is proved that  for the upper deficiency number
$n_+(=n_-)$ of operator $L_0$ the formula $n_+=2m+p$ is true
if $p \in \{1,2, \ldots, m \}$. Moreover, in~\cite{Oro5} a hypothesis about the equality $n_+=4m-p$ for $p \in \{m+1,m+2, \ldots, 2m-1 \}$ is formulated. In~\cite{DolMirz} this hypothesis is proven.

Y.B. Orochko also studied the problem of determining the deficiency numbers of the minimal symmetric differential operator generated by the irregular differential expression~(\ref {osn}) on $I$, whose coefficient $c(x)$ has the different  orders of zero  $p>0$, $q>0$ where $p$, $q$ are not integers. For correct definition of the operator an additional condition is required, namely, it is necessary and sufficient that $\min\{p,q\}> m- \frac 1 2$.
In his works, Y.B.Orochko obtained some estimates for the deficiency numbers of the operator in this case.
 Examining this problem, he only supposed that the coefficient $c(x)$ is differentiated a sufficient number of times and relied on the well-developed asymptotic methods of the theory of ordinary linear differential equations; however, because of the specificity of the differential expression, the complexity of calculations when determining the fundamental system of solutions of the corresponding equations increases considerably as  the order of the equation increases.
In this paper, this fact will also form the basis of the assumption that $a(z)$ and $b(z)$ are analytic functions, as in that case, it is already possible to use methods of the analytical theory of differential equations. Such approach was also used in \cite{DolMirz} and \cite{Braeu}.

\section{Auxiliary results}
Let us formulate some basic facts required in the following sections.

\begin{lemma}  For any function $f(x) \in D_0$ the following is true: \\
$A.$ if $p \in \{1,2, \ldots ,m \}$, then \\
$1.$ $f^{[l]}(0)=0$ and $f^{[l]}(x)=O(x^{p+m-l})$ as $x \to +0$
if $l=m, m+1, \ldots, p+m-1$\\
$2.$ for any $2m-p$ complex numbers $b_l$, $l=0,1,
\ldots, m-1$ and $l=m+p,m+p+1, \ldots, 2m-1$ one can find a function
$f(x) \in D_0$ so that $f^{[l]}(0)=b_l$ for all these values of $l$\\
$B.$ if $p \in \{m+1,m+2, \ldots ,2m-1 \}$, then \\
$1^{\prime}.$ $f^{[l]}(0)=0$ and $ f^{[l]}(x)=O(x^{p+m-l})$ as $x
\to +0$ if $l=m, m+1, \ldots, 2m-1$\\
$2^{\prime}.$ for $m$ arbitrary complex numbers $b_l$,
$l=0, 1, \ldots, m-1$, one can find a
function $f(x) \in
D_0$ such that $f^{[l]}(0)=b_l$ for all specified values of index $l$.
\end{lemma}

The validity of this lemma follows directly from the definition of quasi-derivatives and in case $p \leq m$ is included in~\cite{Oro5}, $p > m$ is in \cite{DolMirz}.

Consider now two auxiliary differential expressions
$$l_{2m,-}[f](x) = (-1)^m ( (-x)^q b(x) f^{(m)} )^{(m)}(x), \;\;x \in [-1,0),$$
$$l_{2m,+}[f](x) = (-1)^m ( x^p a(x) f^{(m)} )^{(m)}(x), \;\; x \in (0,1].$$
Fix a complex number $\lambda$, $\Im \lambda \ne 0$, and consider the differential equations
\begin {equation}
\label {C} l_{2m,-}[y](x) = \lambda y(x), \;\; x \in [-1,0)
\end {equation}
\begin {equation}
\label {B} l_{2m,+}[y](x) = \lambda y(x), \;\; x \in (0,1].
\end {equation}

We denote by $N_{pq}$ the deficiency subspace of the symmetric operator $L_0^{pq}$,
corresponding to $\overline \lambda$, and by $N_{q}$ and $N_{p}$ the deficiency subspaces of the auxiliary
symmetric minimal operators $L_{0}^{q}$ and $L_{0}^{p}$, which were generated by the differential expressions $l_{2m,-}[f](x)$ in $L_2[-1,0)$
and $l_{2m,+}[f](x)$ in $L_2(0,1]$, respectively, and corresponding
to the same $\overline \lambda$.
The deficiency numbers of the operators generated by $l_{2m,-}$ and
$l_{2m,+}$ - $dim \; N_{q}$ and $dim \; N_{p}$ are determined by the maximum number of the
linearly independent solutions of the equations~(\ref {C}) and~(\ref {B}) in the space $L_2[-1,0)$ and $L_2(0,1]$,  respectively.

Let $y_+(x)$ and $y_-(x)$ be restrictions of the function $y(x)$ defined
for $x \in I$ to $(0,1]$ and $[-1,0)$ respectively.

The following lemma is true (see~\cite{Oro5}).
\begin{lemma} \label{lemma2}
For any positive integers $p$ and $q$ the deficiency subspace $N_{pq}$ of the
operator $L_{0}^{pq}$ is the lineal of functions $y(x)\in L_2[-1,1]$ having the three properties:
\begin{enumerate}
\item $y_{-}(x) \in N_{q}$, $y_{+}(x) \in N_{p}$,
\item for each function $f(x) \in D_0$ there exist one-sided limits $[f,y](-0)$ and $[f,y](+0)$,
\item the conjugation condition $[f,y](-0) = [f,y](+0)$ is fulfilled at $x = 0$ for each $f(x) \in D_0$.
\end{enumerate}
\end{lemma}

This lemma implies that in order to find the deficiency numbers of the operator $L_{0}^{pq}$ it is necessary to calculate the limits of sesquilinear forms $[f,y](-0)$ and $[f,y](+0)$, and for this purpose it is enough to define the limits $[f,y_{k,-}](-0)$ and $[f,y_{k,+}](+0)$, where functions $ y_{k,-}$ form the basis of $N_{q}$ and $y_{k,+}$ - the basis of $N_{p}$.

Therefore we must define the bases of spaces $N_{q}$ and $N_{p}$ with such precision that it becomes possible to calculate the above mentioned limits. We note that  quasi-derivatives of the functions $y_{k,-}$ и $y_{k,+}$ are included in the sesquilinear forms $[f,y_{k,-}]$ and $[f,y_{k,+}]$ and in order to define their behavior in the neighborhood of zero  we need to know some of the terms of the asymptotic solutions of the equations~(\ref {C}) and~(\ref {B}) and their quasi-derivatives. On the other hand,  the assumption made about analyticity of the functions $a(x)$ and $b(x)$ allows the construction of exact solutions of the corresponding equations, namely, under our assumptions, equations~(\ref {C}) and~(\ref {B}) are equations with only one regular singular point $x=0$, so by applying the Frobenius' method
(see ~\cite[Ch. XVI]{Ince}), we can construct a fundamental system of solutions of this equation. This will be done in Lemmas~\ref {lemma3} and~\ref {lemma4}.

 It should be noted that the process of constructing a fundamental system essentially depends on the order of the zero of the coefficient of the equation, i.e. on number $p$. For this reason, it is useful to divide values of $p$ into 2 groups, namely, $p \in \{1,2 \ldots m \}$ and $p \in \{m+1,m+2 \ldots 2m-1 \}$.

The following lemmas are true:
\begin{lemma} \label{lemma3}
  For $p=1, 2, \ldots, m$ the differential equation~(\ref {B}) has a fundamental system of solutions $y_{0,+}, y_{1,+}, \ldots, y_{2m-1,+} $, so that all solutions of this system belong to the space $L_2(0,1]$ and are determined by the formulae :
$$y_{i,+}=x^{2m-p-i-1} \left ( \sum\limits_{\nu=0}^{\infty}\gamma_{\nu}x^{\nu} \right ), \quad (0 \le i \le m-p-1),$$
$$y_{m-p+2i-2,+}=x^{m-i} \left ( \sum\limits_{\nu=0}^{\infty}\beta_{\nu}x^{\nu} \right ), \quad (1 \le i \le p), $$
$$y_{m-p+2i-1,+} = x^{m-i} \left ( \sum\limits_{\nu=0}^{\infty}\alpha_{\nu}^{1}x^{\nu}+\ln(x)\sum\limits_{\nu=0}^{\infty}\alpha_{\nu}^{2}x^{\nu} \right ), \quad (1 \le i \le p),$$
$$y_{m+p+i,+}=x^{m-p-i-1} \left ( \sum\limits_{\nu=0}^{\infty}\delta_{\nu}^1 x^{\nu} + \ln(x) \sum\limits_{\nu=i+1}^{\infty}\delta_{\nu}^2x^{\nu} \right ), \quad (0 \le i \le m-p-1),$$
\noindent where $\alpha_{\nu}^{1}, \alpha_{\nu}^{2},\beta_{\nu},\gamma_{\nu},\delta_{\nu}^1, \delta_{\nu}^2$ are numbers depending on $m,p,i$ and the serial expansion coefficients of the function $a(x)$.
\end{lemma}
\begin{lemma} \label{lemma4}
  For $p=m+1, m+2, \ldots, 2m-1$ the differential equation~(\ref {B}) has a fundamental system of solutions $y_{0,+}, y_{1,+}, \ldots, y_{2m-1,+} $ so that
  $3m-p$ of the solutions belong to $L_2(0,1]$ and are determined by formulae:
 $$y_{i,+} = x^{m-1-i} \left ( \sum\limits_{\nu=0}^{\infty}\gamma_{\nu}x^{\nu} \right ), \quad (0 \le i \le p-m-1),$$
 $$y_{3m-p-2i-2,+} = x^i \left ( \sum\limits_{\nu=0}^{\infty}\alpha_{\nu}^{0}x^{\nu}+\sum\limits_{\nu=2m-p}^{\infty}\alpha_{\nu}^{1}x^{\nu} \ln x+ \sum\limits_{\nu=2m-p+1}^{\infty}\alpha_{\nu}^{3}x^{\nu} \ln^2x+ \ldots + \right.$$
 $$\left. \sum\limits_{\nu=m-i}^{\infty}\alpha_{\nu}^{p-m-i-1}x^{\nu} \ln^{p-m-i-1}x \right ) , \quad (0 \le i \le 2m-p-1), $$
 $$y_{3m-p-2i-1,+} = x^i \left ( \sum\limits_{\nu=0}^{\infty}\beta_{\nu}^{0}x^{\nu}+\sum\limits_{\nu=0}^{\infty}\beta_{\nu}^{1}x^{\nu} \ln x+\sum\limits_{\nu=2m-p}^{\infty}\beta_{\nu}^{2}x^{\nu}\ln^2x+ \ldots + \right. $$
 $$ \left. \sum\limits_{\nu=m-i}^{\infty}\beta_{\nu}^{p-m-i}x^{\nu} \ln^{p-m-i}x \right ), \quad (0 \le i \le 2m-p-1),$$

 \noindent where $\alpha_{\nu}^{j}, \beta_{\nu}^{j},\gamma_{\nu} $ are numbers which are dependent on $m,p,i$ and the serial expansion coefficients of the function $a(x)$.
\end{lemma}

Using the results of Lemmas~\ref {lemma3} and~\ref {lemma4}, we can prove the following lemma.

\begin{lemma} \label{lemma5}
Let $f(x) \in D_0$. Then \\
 $1)$ for $p \in \{1, 2, \ldots, m \}$ equation~(\ref {B}) has a fundamental system of solutions $y_{k,+}$, $k=1,2, \ldots, 2m$ whose elements belong to the space $L_2(0,1]$ and for which values $[f,y_{k,+}](+0)$ can be calculated by the formula
$$
  [f,y_{k,+}](+0)=
  \begin{cases}
\alpha_k f^{[2m-k]}(0), & \text{ if $k=1,2, \ldots, m-p$,}\\
0, & \text{ if $k=m-p+1,m-p+2 \ldots, m$,}\\
\alpha_k f^{[2m-k]}(0), & \text{ if $k=m+1,m+2, \ldots, 2m$,}
 \end{cases}
$$
$2)$ for $p \in \{m+1,m+2, \ldots, 2m-1 \}$ equation~(\ref {B}) has a fundamental system of solutions $y_{k,+}$, $k=1,2, \ldots, 2m$ whose
 $3m-p$ elements $y_{1,+}, y_{2,+}, \ldots, y_{m,+},$ $y_{p+1,+},y_{p+2,+}, \ldots, y_{2m,+}$ belong to $L_2(0,1]$ and for which the values $[f,y_{k,+}](+0)$ can be calculated by the formulae
$$
  [f,y_{k,+}](+0)=
  \begin{cases}
\beta_{k} f^{[2m-k]}(0), & \text{ if $k=p+1,p+2, \ldots, 2m$,}\\
0, & \text{ in all other cases.}
 \end{cases}
$$
Here  $\alpha_k$, $k \in \{1,2, \ldots, 2m \} \backslash
\{m-p+1,m-p+2, \ldots, m \}$ and $\beta_{k}$, $k=p+1,p+2,
\ldots, 2m$ are non-zero constants.
\end{lemma}
\noindent \textbf {Remark 2.1.} We point out that the substitution $x \rightarrow -x$ reduces equation~(\ref {C}) to an equation of the form~(\ref {B}), therefore in the space
$N_{q}$ there is a basis $y_{k,-} (x)$  whose elements have
the properties listed in Lemma~\ref {lemma5}.

\noindent \textbf {Remark 2.2.} We also mention here that by Lemmas \ref{lemma3} and \ref{lemma4} the solutions $y_{k,+}$ and $y_{k,-}$ are entire in $\lambda$ and the main terms of their asymptotic and quasi-derivatives  do not depend on $\lambda$. Hence we can assume $\lambda =0$ in some situations below.

\section{Deficiency numbers of the operator $L_{0}^{pq}$}
Let us now formulate and prove the main theorem about the deficiency numbers of the operator $L_0^{pq}$.

\begin{theorem}
\label{maintheor}
The deficiency numbers of the operator $L_0^{pq}$ are defined by the formula:
$$
n_{pq}= \begin{cases}
4m-\max\{p,q\}, & \text {if $p,q \in \{m+1,m+2, \ldots, 2m-1\}$, } \\
2m+\min\{p,q\}, & \text {if $p,q \in \{1,2, \ldots, m\}$, } \\
3m+p-q, & \text {if $p \in \{1,2, \ldots, m \}$, $ q \in \{m+1,m+2, \ldots, 2m-1\}$. }
\end{cases}
$$
\end{theorem}
\begin{proof}
The proof scheme of this theorem is the same for all cases and restates the arguments presented in ~\cite{Oro5} for the case $p=q=1,2, \ldots m$.

Above we have determined the structure of linearly independent solutions of the  differential equations~(\ref {B}) and~(\ref {C}) belonging to $L_2(0,1]$ and $L_2[-1,0)$ respectively (Lemmas~\ref {lemma3} and ~\ref {lemma4}), and also determined the values of the limits of the sesquilinear forms corresponding to these solutions (Lemma~\ref {lemma5}).

Let us assume that $p,q \in \{m+1,m+2, \ldots, 2m-1\}$.
A set of functions
from $L_2[-1,1]$ with property
1 of Lemma~\ref {lemma2} is a class of functions $y(x)$ defined on $[-1,1]$ for which the following representations are true:
$$
 y_{-}(x) = \sum\limits_{k=1}^{m}d_k y_{k,-}
(x)+\sum\limits_{k=q+1}^{2m}d_k y_{k,-} (x), \;\;x \in
[-1,0),
$$
\begin {equation}
\label {D} y_{+}(x) = \sum\limits_{k=1}^{m}c_k y_{k,+}
(x)+\sum\limits_{k=p+1}^{2m}c_k y_{k,+} (x), \;\;x \in (0,1],
\end {equation}
 for some complex coefficients $c_k$ $(k=1,2, \ldots, m,p+1,p+2, \ldots ,
2m)$, $d_k$ $(k=1,2, \ldots, m,q+1,q+2, \ldots ,
2m)$, where $y_{k,-} (x)$ and $y_{k,+}(x)$ are the  bases of the spaces
$N_{q}$ and $N_{p}$ with the properties listed in Lemma~\ref {lemma5} .

  By (\ref {D}) and the linearity of the form $[f,g](x)$, it follows that for
each function  $y(x)$ and any function $f(x) \in D_0$ we have
      $$[f,y_{-}](-0) = \sum\limits_{k=q+1}^{2m}d_{k} [f,y_{k,-}](-0),$$
    $$[f,y_{+}](+0) =\sum\limits_{k=p+1}^{2m}c_{k} [f,y_{k,+}](+0).$$
 By Lemma~\ref {lemma5}, we determine that
    $$[f,y_{-}](-0) = \sum\limits_{k=q+1}^{2m}d_{k} \beta_{k} f^{[2m-k]}(0),$$
    $$[f,y_{+}](+0) =\sum\limits_{k=p+1}^{2m}c_{k} \alpha_{k} f^{[2m-k]}(0),$$
where $\beta_{k}$ and $\alpha_{k}$ are non-zero numbers.

 Let us assume for definiteness that $q>p$.
 From the set of functions $y(x)$ we now select the class of those which additionally possess property 3 of Lemma~\ref {lemma2}. This fulfillment of this condition for $y(x)$ with any function $f(x) \in D_0$ is equivalent to the following system
\begin {gather}
 \alpha_{k}c_{k}=\beta_{k}d_{k}, \;\; k=q+1,q+2,
 \ldots, 2m, \notag \\
 \label {G} \alpha_{k}c_{k}=0, \;\; k=p+1,p+2,
 \ldots, q,
\end {gather}
on $6m-p-q$ coefficients $c_k$ $(k=1,2, \ldots, m,p+1,p+2, \ldots, 2m)$ and $d_k$
 $(k=1,2, \ldots, m,q+1,q+2, \ldots, 2m)$.

Hence, by Lemma~\ref {lemma2} the deficiency subspace $N_{pq}$ of the operator $L_{0}^{pq}$
 is the lineal of functions $y(x)$ admitting the representation~(\ref {D}) with the coefficients $c_k$,
$(k=1,2, \ldots,m,p+1,p+2, \ldots, 2m)$ and $d_k$
 $(k=1,2, \ldots, m,q+1,q+2, \ldots, 2m)$
which are connected by relations~(\ref {G}), but otherwise arbitrary.
 The dimension of this subspace is equal to the number of those listed $6m-p-q$ coefficients that can assume arbitrary values under these constraints.

 In this relation~(\ref {G}) \\
$-$ there are no coefficients $c_{k},d_{k}$ for $k=1,2, \ldots, m$, consequently these $2m$ coefficients are arbitrary,\\
$-$ the coefficients $c_{k}=0$ for $k=p+1,p+2, \ldots, q$,\\
$-$ among the coefficients $c_{k},d_{k}$ for $k=q+1,q+2, \ldots, 2m$
only $2m-q$ coefficients are arbitrary.

Hence, the total number of the coefficients taking arbitrary values equals
 $2m+2m-q=4m-q$, so, in this case, the deficiency number $n_{pq}$
 of the operator $L_0^{pq}$ is defined by the formula $n_{pq}=4m-q$.

In the case $q<p$, to repeat the arguments above, we have the equality
$n_{pq}=4m-p$. Therefore, $n_{pq}=
4m-\max\{p,q\}.$

The cases $p,q \in \{1,2, \ldots, m\}$ and $p \in \{1,2, \ldots, m\}$, $ q \in \{m+1,m+2, \ldots, 2m-1\}$ can be proven similarly.
\end{proof}

\section{Self-adjoint extensions of the operator $L_0^{pq}$}
It is well known that the classification of the self-adjoint extensions of $L_0^{pq}$ depends, in an essential way, on the deficiency numbers of $L_0^{pq}$.

At first we summarize a few properties of the basis of the deficiency space $N_{pq}$ in the form of lemma.
\begin{lemma}
\label{lemma7}
 1. If $p,q \in \{m+1,m+2, \ldots, 2m-1\}$, then the basis of $N_{pq}$, corresponding to $\overline{\lambda}$ consists of $4m-\max\{p,q\}$ functions
$$
\phi_j =
 \begin{cases}
c_{j,-} z_{j,-}, \quad x \in [-1,0), \\
c_{j,+} z_{j,+}, \quad x \in (0,1],
 \end{cases}
$$
where $z_{j,+}=y_{k,+}, z_{j,-}=0$ if $j(=k)=1,2, \ldots, m$ and $z_{j,+}=0, z_{j,-}=y_{k,-}$ if $j=m+1,m+2, \ldots, 2m, k=1,2, \ldots, m$ and
$z_{j,+}=y_{k,+}, z_{j,-}=y_{k,-}$ if $j=2m+1,m+2, \ldots, 4m-\max\{p,q\}, k=\max\{p,q\}+1,\max\{p,q\}+2, \ldots, 2m$, $c_{j,-}, c_{j,+}$ are real numbers and the conjugation condition
$[f,z_{j,-}](-0) = [f,z_{j,+}](+0)$ is fulfilled at $x = 0$ for each $f(x) \in D_0$ if $j=1,2, \ldots, 4m-\max\{p,q\}$.

2.
If $p,q \in \{1,2, \ldots, m\}$, then the basis of $N_{pq}$, corresponding to $\overline{\lambda}$ consists of $2m+\min\{p,q\}$ functions
$$
\phi_j =
 \begin{cases}
c_{j,-} z_{j,-}, \quad x \in [-1,0), \\
c_{j,+} z_{j,+}, \quad x \in (0,1],
 \end{cases}
$$
 where $z_{j,+}=y_{k,+}, z_{j,-}=y_{k,-}$ if $j(=k)=1,2, \ldots, m-\max\{p,q\}$ and  $z_{j,+}=y_{k,+}, z_{j,-}=0$ if $j=m-\max\{p,q\}+1,m-\max\{p,q\}+2, \ldots, m, k=m-\max\{p,q\}+1,m-\max\{p,q\}+2, \ldots, m$ and $z_{j,+}=0, z_{j,-}=y_{k,-}$ if $j=m+1,m+2, \ldots, m+\min\{p,q\}, k=m-\min\{p,q\}+1,m-\min\{p,q\}+2, \ldots, m$ and
$z_{j,+}=y_{k,+}, z_{j,-}=y_{k,-}$ if $j=m+\min\{p,q\}+1,m+\min\{p,q\}+2, \ldots, 2m+\min\{p,q\}, k=m+1,m+2, \ldots, 2m$, $c_{j,-}, c_{j,+}$ are real numbers and the conjugation condition $[f,z_{j,-}](-0) = [f,z_{j,+}](+0)$ is fulfilled at $x = 0$ for each $f(x) \in D_0$ if $j=1,2, \ldots, 2m+\min\{p,q\}$.

3. If $p \in \{1,2, \ldots, m\}$, $ q \in \{m+1,m+2, \ldots, 2m-1\}$, then the basis of $N_{pq}$, corresponding to $\overline{\lambda}$ consists of $3m+p-q$ functions
$$
\phi_j =
 \begin{cases}
c_{j,-} z_{j,-}, \quad x \in [-1,0), \\
c_{j,+} z_{j,+}, \quad x \in (0,1],
 \end{cases}
$$
where $z_{j,+}=0, z_{j,-}=y_{k,-}$ if $j(=k)=1,2, \ldots, m$ and $z_{j,+}=y_{k,+}, z_{j,-}=0$ if $j=m+1,m+2, \ldots, m+p, k=m-p+1,m-p+2, \ldots, m$ and
$z_{j,+}=y_{k,+}, z_{j,-}=y_{k,-}$ if $j=m+p+1,m+p+2, \ldots, 3m+p-q, k=q+1,q+2, \ldots, 2m$, при $c_{j,-}, c_{j,+}$ are real numbers and the conjugation condition $[f,z_{j,-}](-0) = [f,z_{j,+}](+0)$ is fulfilled at $x = 0$ for each $f(x) \in D_0$ if $j=1,2, \ldots, 2m+\min\{p,q\}$.
\end{lemma}
The proof of this lemma directly follows from Lemmas \ref{lemma3} and \ref{lemma4} and the proof of Theorem 3.1

\noindent \textbf {Remark 4.1} A similar result holds for the deficiency space $N_{pq}$, corresponding to $\lambda$.

Now, we can use the Gram--Schmidt orthogonalization process to $\phi_j$ and construct an orthonormal basis for  $N_{pq}$ corresponding to $\overline{\lambda}$ and apply Theorem 2 (see  \cite[\S 18]{Naim}) to our case.

Let
$$
\varphi_1(x),\varphi_2(x), \ldots, \varphi_n(x)
$$
be an orthonormal basis in $N_{pq}$, corresponding ${\overline {\lambda}}$. Then the functions
$$
-\overline{\varphi_1(x)},-\overline{\varphi_2(x)}, \ldots, -\overline{\varphi_n(x)}
$$
form an orthonormal basis in $N_{pq}$, corresponding $\lambda$.

\begin{theorem}
\label{T1}
Every self-adjoint extension $L_u^{pq}$ of the operator $L_0^{pq}$ with the deficiency numbers $(n,n)$ can be characterized by means of a unitary $n \times n$ matrix $u=[u_{\mu \nu}]$ in the following way. \\
The domain of definition $D_{u}$ is the set of all functions $z(x)$ of the form
$$
z(x)=y(x)+\psi(x),
$$
where $y(x) \in D_0$,  $\psi(x)$ is a linear combination of the functions
$$
\psi_{\mu}(x)=\varphi_{\mu}(x)+\sum\limits_{\nu=1}^n u_{\nu \mu} \overline{\varphi_{\nu}(x)}, \\ \mu=1,2, \ldots, n
$$
and
$$
n=
 \begin{cases}
4m-\max\{p,q\}, & \text{ если $p,q \in \{m+1,m+2, \ldots, 2m-1\}$;}\\
2m+\min\{p,q\}, & \text{ если $p,q \in \{1,2, \ldots, m\}$;}\\
3m+p-q,   & \text{ если $p \in \{1,2, \ldots, m\}$, $ q \in \{m+1,m+2, \ldots, 2m-1\}$.}
 \end{cases}
$$
\noindent Conversely, every unitary $n \times n$ matrix $u=[u_{\mu \nu}]$ determines in the way described above a certain self-adjoint extension $L_u^{pq}$ of the operator $L_0^{pq}$.
\end{theorem}
The following theorem characterizes the domain of definition of the operator $L_u^{pq}$ by means of boundary conditions.
\begin{theorem}
\label{SD}
The domain of definition $D_{u}$ of an arbitrary self-adjoint extension $L_u^{pq}$ of the operator $L_0^{pq}$ with the deficiency indices $(n,n)$ consists of the set of all functions $y(x) \in D$, which satisfy the conditions
\begin{equation}
\label{BC1}
[y,w_k](1)-[y,w_k](+0)+[y,w_k](-0)-[y,w_k](-1)=0, k=1,2, \ldots, n,
\end{equation}
where $w_1,w_2, \ldots, w_n$ are certain functions belonging to $D$ and determined by $L_u^{pq}$, which are linearly independent modulo $D_0$ and for which the relations
\begin{equation}
\label{BC}
[w_j,w_k](1)-[w_j,w_k](+0)+[w_j,w_k](-0)-[w_j,w_k](-1)=0, k=1,2, \ldots, n
\end{equation}
hold.
Conversely, for arbitrary functions $w_1,w_2, \ldots, w_n$ belonging to $D$ which are linearly independent modulo $D_0$ and which satisfy the relations (\ref{BC}), the set of all functions $y(x) \in D$ which satisfy the conditions (\ref{BC1}) is the domain of definition of a self-adjoint extension of the operator $L_0^{pq}$.
\end{theorem}
The proof of these two theorems is exactly repeat the proof of Theorems 2 and 4 in \cite [\S 18, \S 18.1] {Naim}.

Next, using the ideas of \cite{Fult} and the approach of \cite{Sun}, \cite{Zettl2} and \cite{Sun1} we can specialize the domain of definition $D_{u}$ of an arbitrary self-adjoint extension $L_u^{pq}$ in more applicable way.

Let $n_{pq}$ denote the deficiency index of $L_0^{pq}$ as above.
Let $d_1$ and $d_2$ be numbers of linearly independent solutions of (\ref{C}) and (\ref{B}) which form the basis of the deficiency subspace $N_{pq}$ and $d_3$ be a number of these functions satisfying the conjugation condition. We note here that $d_1+d_2-d_3=n_{pq}$.

Assume that the functions $ \phi_1(x,\overline{\lambda}), \phi_2(x,\overline{\lambda}), \ldots, \phi_{n_{pq}}(x,\overline{\lambda})$ form the basis of $N_{pq}$, corresponding $\overline{\lambda}$  and $ \psi_1(x,\lambda), \psi_2(x,\lambda), \ldots, \psi_{n_{pq}}(x,\lambda)$ form the basis of $N_{pq}$, corresponding $\lambda.$
We note here that
$$
\label{fphi}
\begin{aligned}
& \phi_i(x,\overline{\lambda}) =
 \begin{cases}
\phi_{j,-}(x,\overline{\lambda}), \quad x \in [-1,0) ,  \qquad i,j=1,2, \ldots, d_1 - d_3,\\
0, \quad x \in (0,1],
 \end{cases}
\\
& \phi_i(x,\overline{\lambda}) =
 \begin{cases}
 \phi_{j,-}(x,\overline{\lambda}), \quad x \in [-1,0) , \qquad i,j=d_1-d_3+1, \ldots, d_1,\\
 \phi_{j,+}(x,\overline{\lambda}), \quad x \in (0,1], \qquad j=d_2-d_3+1, \ldots, d_2,
 \end{cases}
 \\
& \phi_i(x,\overline{\lambda}) =
 \begin{cases}
0, \quad x \in [-1,0) , \qquad i=d_1+1 , \ldots, n_{pq}, \\
\phi_{j,+}(x,\overline{\lambda}), \quad x \in (0,1], \qquad j=1,2 \ldots d_2-d_3,
 \end{cases}
 \end{aligned}
$$
where $\phi_{j,-}(x,\overline{\lambda})$ $(j=1,2, \ldots, d_1)$ and $\phi_{j,+}(x,\overline{\lambda})$ $(j=1,2, \ldots, d_2)$ are the linearly independent solutions of (\ref{C}) and (\ref{B}) which form the basis of the deficiency subspace $N_{pq}$ corresponding $\overline{\lambda}$.
We have a similar representation for the functions $\psi_i(x,\lambda)$, $(i=1,2, \ldots, n_{pq})$

For the convenience, we denote
\begin{equation}
\label{bas}
\chi_{2i-1}=\psi_i(x,\lambda), \chi_{2i}=\phi_i(x,\overline{\lambda}), i=1,2 \ldots n_{pq}
\end{equation}
and $\chi_{j,-}$, $\chi_{j,+}$ are restrictions of the function $\chi_{j}$ $(j=1,2, \ldots, 2n_{pq})$ to $[-1,0)$  and $(0,1]$ respectively.

Let $g_{i,-}(x)$ and $g_{i,+}(x)$ $(i=1, \ldots,2m)$ be sets of functions in $D$ defined on $[-1,1]$, which satisfy the following conditions:
\begin{equation}
\label{g_func}
\begin{aligned}
& g_{i,-}^{[k-1]}(-1)=\delta_{ik}, \;\; g_{i,-}^{[k-1]}(a)=0,\;(-1<a<0), \; i,k=1, \ldots, 2m, \\
& g_{i,-}(x)=0, \; x \ge a, \\
& g_{i,+}^{[k-1]}(1)=\delta_{ik}, \;\; g_{i,+}^{[k-1]}(b)=0,\;(0<b<1), \; i,k=1, \ldots, 2m, \\
& g_{i,+}(x)=0, \; x \le b. \\
\end{aligned}
\end{equation}
By the Naimark Patching Lemma (\cite[Chap.5,\S 17]{Naim}), there exist such functions.

Since  $g_{i,-}(x) \in D$ $(i=1, \ldots,2m)$, by Theorem \ref{T1}, we have
\begin{equation}
\label{g_func1}
g_{i,-}=y_{0i,-} + \sum\limits_{j=1}^{2n_{pq}} a_{ij,-} \chi_{j}, \;\; y_{0i,-} \in D_0 \; (i=1,\ldots,2m).
\end{equation}\\
Also, we have
\begin{equation}
\label{g_func2}
g_{i,+}=y_{0i,+} + \sum\limits_{j=1}^{2n_{pq}} a_{ij,+} \chi_{j}, \;\; y_{0i,+} \in D_0 \; (i=1,\ldots,2m).
\end{equation}
Similarly, as it was done in Lemma 1 in \cite{Sun1} we can show that  $$rank \; X_- (:=([\chi_i,\chi_j](-0))_{2n_{pq} \times 2n_{pq}})=2d_1-2m$$ and $$rank \; X_+ (:=([\chi_i,\chi_j](+0))_{2n_{pq} \times 2n_{pq}})=2d_2-2m.$$
Therefore, it is possible to arrange the functions $\chi_i$, so that the matrices $X_-$ and $X_+$ can be represented as
\begin{equation}
\label{MatrX}
X_-:=
\left(\begin{array} {l}
X^{1,-}_{2d_1 -2m \times 2n_{pq}} \\
X^{2,-}_{2m \times 2n_{pq}} \\
X^{3,-}_{2d_2-2d_3 \times 2n_{pq}}
\end{array} \right), \qquad
X_+:=
\left( \begin{array}{l}
X^{1,+}_{2d_2-2m \times 2n_{pq}} \\
X^{2,+}_{2m \times 2n_{pq}} \\
X^{3,+}_{2d_1-2d_3 \times 2n_{pq}}
\end{array} \right ),
\end{equation}
where $rank \; X^{1,-}= 2d_1 -2m$ and $rank \; X^{1,+} =2d_2-2m$. Let $$A_-(:=(a_{ij,-})_{2m \times 2n_{pq}})=(C^-_{2m \times (2d_1-2m)}D^-_{2m \times 2m}F^-_{2m \times (2d_2-2d_3)})$$ and $$A_+(:=(a_{ij,+})_{2m \times 2n_{pq}})=(C^+_{2m \times (2d_2-2m)}D^+_{2m \times 2m}F^+_{2m \times (2d_1-2d_3)}).$$ Then using ideas of Lemma 2 in \cite{Sun1} it is easy enough to obtain that $rank \; D^-= 2m$ and $rank \; D^+ =2m$. Therefore the following lemmas take place.
\begin{lemma}
Let $n_1=2d_1-2m$ and $n_2=2d_2-2m$. Suppose $\{ \chi_i \}$ are the functions defined in (\ref{bas}), which satisfy (\ref{MatrX}), then each of the functions $\chi_{i,-}$ $(i=2d_1-2m+1, \ldots 2d_1)$ and $\chi_{i,+}$ $(i=2d_2-2m+1, \ldots 2d_2)$ has a unique representation
\begin{equation}
\label{phi_m}
\begin{array}{l}
\chi_{i,-}=\tilde{y}_{i0,-}+\sum\limits_{j=1}^{2m} c_{j,-} g_{j,-} + \sum\limits_{s=1}^{n_1} b_{is,-} \chi_{s,-}, \\
\chi_{i,+}=\tilde{y}_{i0,+}+\sum\limits_{j=1}^{2m} c_{j,+} g_{j,+} + \sum\limits_{s=1}^{n_2} b_{is,+} \chi_{s,+},
\end{array}
\end{equation}
where $\tilde{y}_{i0,-}, \tilde{y}_{i0,+} \in D_0$ and $g_{j,-}, g_{j,+}$ satisfy (\ref{g_func1}) and (\ref{g_func2}) respectively.
\end{lemma}
\begin{lemma}
\label{Dom} The domain $D$ of the maximal operator $L$ can be represented as
$$
D=D_0+span\{g_{1,-},\ldots, g_{2m,-}\} +span\{\chi_{1,-},\ldots, \chi_{n_1,-}\} +
$$
$$ span\{g_{1,+},\ldots, g_{2m,+}\} +span\{\chi_{1,+},\ldots, \chi_{n_2,+}\}.$$
\end{lemma}
We mention that the method of proof given in \cite{Sun1} (see Theorem 1) can be also adapted to prove Lemma 4.5.

In \cite{Zettl1} it has been shown that in the one singular end-point case the complex-valued functions $\chi_{k,-}$ and $\chi_{k,+}$ in (\ref{phi_m}) can be replaced by the real-valued functions.

Let $E_k=((-1)^r \delta_{r,k+1-s})_{r,s=1}^k$ be a symplectic matrix of the order $k$ and
$u_{i,-}$ $(i=1,2\ldots, n_q)$ and $u_{i,+}$ $(i=1,2,\ldots,n_p)$ be linearly independent solutions of $l_{2m,-}[u](x)=0$ and $l_{2m,+}[u](x)=0$ which lie in $L_2[-1,0)$ and $L_2(0,1]$ respectively.

The solutions $u_{i,-}, i=1, 2 \ldots, n_q$ on $[-1,0)$ can be ordered such that the $n_1 \times n_1$ matrix $U_-=([u_{i,-},u_{j,-}](-1))$, $i,j=1,2, \ldots, n_1$, is given by
$$
U_-=(-1)^{m+1}E_{n_1}
$$
and the solutions $u_{i,+}, i=1, 2 \ldots, n_p$ on $(0,1]$ can be ordered such that the $n_2 \times n_2$ matrix $U_+=([u_{i,+},u_{j,+}](1)), i,j=1,2 ,\ldots, n_2$, is given by
$$
U_+=(-1)^{m+1}E_{n_2}.
$$
Let us determine functions $g_{j,-} \in D$, $j=1, \ldots, 2m$ such that $g_{j,-}(t)=0$ for $t \ge a_- (-1 < a_- < 0)$ and the $2m \times 2m$ matrix  $G_-=([g_{i,-},g_{j,-}](-1)), i,j=1,2 ,\ldots, 2m,$ is given by
$$
G_-=E_{2m}
$$
and
functions $g_{j,+} \in D$, $j=1, \ldots, 2m$ such that $g_{j,+}(t)=0$ for $t \le a_+ (0 < a_+ < 1)$ and the $2m \times 2m$ matrix  $G_+=([g_{i,+},g_{j,+}](1)), i,j=1,2 ,\ldots, 2m,$ is given by
$$
G_+=E_{2m}.
$$

Using the approach of \cite{Zettl1} and Remark 2.2 we have
\begin{lemma}
\label{cor}
Let the numbers $n_1$, $n_2$ and the functions $u_{i,-}$, $u_{i,+}$, $g_{k,-}$, $g_{k,+}$ are determined as above then each $y\in D$ can be uniquely written as
$$
y=y_0 + \sum\limits_{j=1}^{2m} d_{j,-} g_{j,-} + \sum\limits_{k=1}^{n_1} h_{k,-} u_{k,-} + \sum\limits_{j=1}^{2m} d_{j,+} g_{j,+} + \sum\limits_{k=1}^{n_2} h_{k,+} u_{k,+},
$$
where $y_0 \in D_0$, $d_{j,-}$, $h_{k,-}$, $d_{j,+}$, $h_{k,+}$ are the complex numbers
and
$$D = D_0+span\{g_{1,-},\ldots, g_{2m,-}\} +span\{u_{1,-},\ldots, u_{n_1,-}\} +$$
$$ span\{g_{1,+},\ldots, g_{2m,+}\} + span\{u_{1,+},\ldots, u_{n_2,+}\}.$$
\end{lemma}
Based on Lemma \ref{cor} we can give a characterization of all self-adjoint domains in terms of real-valued solutions of $l_{2m,-}[y]=0$ and $l_{2m,+}[y]=0$.
\begin{theorem}
 \label{SAD}
 Let the complex numbers $n_1$ and $n_2$ as defined above. A linear submanifold $D_{u}$ of $D$ is the domain of a self-adjoint extension $L_u^{pq}$ of $L_0^{pq}$
if and only if there exists complex $n_{pq} \times 2m$ matrices $A_1$ and $A_2$, a complex $n_{pq} \times n_1$ matrix $B_1$ and a complex $n_{pq} \times n_2$ matrix $B_2$ such that the following conditions hold: \\
1. $rank(A_1,B_1,B_2,A_2) = n_{pq}$; \\
2. $A_2E_{2m}A_2^*-B_2E_{n_1}B_2^*+B_1E_{n_2}B_1^*-A_1E_{2m}A_1^*=0$; \\
3. For each $f \in D_0$
$$
B_1\begin{pmatrix}
[f,u_{1,-}](0-) \\
\vdots \\
[f,u_{n_1,-}](0-)
\end{pmatrix}+
 B_2\begin{pmatrix}
[f,u_{1,+}](0+) \\
\vdots \\
[f,u_{n_2,+}](0+)
\end{pmatrix}=\begin{pmatrix} 0 \\ \vdots \\ 0 \end{pmatrix};
$$
4.  $$A_1 \begin{pmatrix} y_-(-1) \\ \vdots  \\ y_-^{[2m-1]}(-1) \end{pmatrix} + B_1  \begin{pmatrix} [y_-, u_{1,-}](0-)  \\ \vdots \\ [y_-, u_{n_1,-}](0-) \end{pmatrix} + $$
$$+B_2 \begin{pmatrix} [y_+, u_{1,+}](0+)  \\ \vdots \\ [y_+, u_{n_2,+}](0+) \end{pmatrix} + A_2  \begin{pmatrix} y_+(1) \\ \vdots \\ y_+^{[2m-1]}(1)  \end{pmatrix}= \begin{pmatrix} 0 \\ \vdots \\ 0 \end{pmatrix}.$$
\end{theorem}

\section{Spectrum of the self-adjoint extensions of the operator $L_0^{pq}$}
In order to describe the spectrum of each self-adjoint extension we need, in particularly, the following lemma (see~\cite{HinL}).

\begin{lemma} \label{lemma6}
The spectrum of each self-adjoint extension of the operator $L_0$ induced by the differential expression $l_{2n}[f](x) = (-1)^n ( p(x) f^{(n)} )^{(n)}(x), \;\; x \in (0,1]$ is discrete and bounded below if and only if
\begin{equation}
\label{eqspec}
\lim\limits_{x \to 0+} x^{1-2m} \int\limits_0^s s^{4n-2}p(s)^{-1} ds = 0.
\end{equation}
\end{lemma}

Let us now formulate and prove

\begin{theorem}
\label{spect}
The spectrum of any self-adjoint extension of the operator $L_0^{pq}$ is discrete.
\end{theorem}

\begin{proof}
In order to examine the spectrum of self-adjoint extensions of the operator $L_0^{pq}$ we will use the splitting method (see, for example,~\cite{Naim},~\cite{Glaz}).

We will analyze the orthogonal decomposition $L_2[-1,1]=L_2[-1,0) \oplus L_2(0,1]$, where $L_2[-1,0)$, $L_2(0,1]$ are considered as subspaces in $L_2[-1,1],$ consisting of the functions $f(x) \in L_2[-1,1]$ equal to zero respectively if $ x \in (0,1]$ and $x \in [-1,0)$. Let us define the orthogonal sum $L_0^q \oplus L_0^p $ of the operator $L_0^q$ acting in $L_2[-1,0)$ and the operator  $L_0^p$ acting in $L_2(0,1]$, which is a real symmetric operator in  $L_2[-1,1]$. It is obvious that the operator $L_0^{pq}$ is a symmetric extension of the operator $L_0^q \oplus L_0^p $.

Furthermore, we extend the operators  $L_0^q$ and $L_0^p$ into self-adjoint operators $L_{0,u}^q$ and $L_{0,u}^p$ in the spaces $L_2[-1,0)$ and $L_2(0,1]$ respectively, then the direct sum $A = L_{0,u}^q \oplus L_{0,u}^p$ will be a self-adjoint extension of the operator $L_0^q \oplus L_0^p $ and the spectrum of the operator $A$ will be the set-theoretic sum of the spectra of $L_{0,u}^q$ and $L_{0,u}^p$.

On the other hand, the deficiency numbers of the operator $L_0^q \oplus L_0^p $ are finite, and thus, all its self-adjoint extensions have one and the same continuous spectrum. Such extensions are the operator $A$ as well as each self-adjoint extension $L_{0,u}^{pq}$ of the operator $L_0^{pq}$, and therefore, the continuous parts of the spectrum of the two operators $A$ and $L_{0,u}^{pq}$ coincide.

Let us now define the spectrum of operators $L_{0,u}^{q}$ and $L_{0,u}^{p}$. To do this we will use Lemma~\ref {lemma6}.

 We will verify the fulfillment of condition~(\ref {eqspec}) for the coefficient of the differential expression~(\ref{B}). In this case $p(x) = x^p a(x) = x^p \sum\limits_{k=0}^{\infty}a_k x^k$ .

Substituting this expression into the left-hand side of ~(\ref {eqspec}), we find that the equality~(\ref {eqspec}) is valid for $0<p<2m$.
Similarly, we obtain the conditions for $q$, i.e. $0 < q < 2m$.

Therefore, in accordance with Lemma~\ref {lemma6}, the spectrum of the operators $L_{0,u}^{q}$ and  $L_{0,u}^{p}$ is discrete if $0 < q < 2m$ and  $0<p<2m$.

 Consequently, the spectrum of the self-adjoint extensions of the operator $L_{0,u}^q \oplus L_{0,u}^p$ and hence also of the self-adjoint extensions of the operator $L_{0}^{pq}$ for the given values $p$ and $q$ is discrete.
\end{proof}

% ------------------------------------------------------------------------

% ------------------------------------------------------------------------
\noindent Irina Braeutigam \\
Department of Mathematical Analysis, Algebra and Geometry, \\
Institute of Mathematics, Information and Space Technologies, \\
Northern (Arctic) Federal University,
163002, Arkhangelsk, Russia\\
e-mail: irinadolgih@rambler.ru, i.braeutigam@narfu.ru

\end {document}